\documentclass[a4paper]{amsart}
\usepackage{amscd,amssymb}

\theoremstyle{plain}
\newtheorem{thm}{Theorem}[section]

\newtheorem{lem}[thm]{Lemma}
\newtheorem{cor}[thm]{Corollary}
\newtheorem*{thOne}{Theorem\,1}
\newtheorem*{thTwo}{Theorem\,3}
\newtheorem*{crl}{Corollary\,2}
\theoremstyle{remark}

\newcommand{\m}{\phantom{-}}

\newcommand{\N}{\mathbb{N}}
\newcommand{\Z}{\mathbb{Z}}
\newcommand{\R}{\mathbb{R}}
\newcommand{\C}{\mathbb{C}\mkern1mu}
\renewcommand{\H}{\mathbb{H}\mkern1mu}
\newcommand{\RP}{\mathbb{R\mkern1mu P}}
\newcommand{\CP}{\mathbb{C\mkern1mu P}}
\newcommand{\Sph}{\mathbb{S}}

\DeclareMathOperator{\sgn}{sgn}

\newcommand{\ip}{\langle \,\cdot\, , \,\cdot \rangle}
\DeclareMathOperator{\tr}{trace}
\DeclareMathOperator{\id}{id}

\DeclareMathOperator{\codim}{codim}
\DeclareMathOperator{\rank}{rank}
\newcommand{\tp}{^{\mathrm{T}}}
\newcommand{\cplxi}{\mathrm{i}}
\newcommand{\1}{{\mathchoice
{\mathrm 1\mskip-4.2mu\mathrm l}{\mathrm 1\mskip-4.2mu\mathrm l}
{\mathrm 1\mskip-3.9mu\mathrm l}{\mathrm 1\mskip-4.0mu\mathrm l}}}
\newcommand{\bmat}{\left(\begin{smallmatrix}}
\newcommand{\emat}{\end{smallmatrix}\right)}
\newcommand{\abs}[1]{\vert #1\vert}

\newcommand{\SO}{\mathrm{SO}}
\newcommand{\OO}{\mathrm{O}}
\newcommand{\SU}{\mathrm{SU}}
\newcommand{\U}{\mathrm{U}}

\newcommand{\Syp}{\mathrm{Sp}}

\newcommand{\T}{\mathrm{T}}

\DeclareMathOperator{\Sq}{Sq}
\begin{document}
\title[Selfmaps of cohomogeneity one manifolds]{Cohomogeneity one manifolds and selfmaps of nontrivial degree}
\author{Thomas P\smash{\"u}ttmann}
\address{Universit\"at Bonn\\ Mathematisches Institut\\ Beringstr. 1\\
  53115 Bonn\\Germany}
\email{puttmann@math.uni-bonn.de}
\thanks{supported by a DFG Heisenberg scholarship and the DFG priority program SPP 1154}

\subjclass[2000]{Primary 57S15, Secondary 55M25}

\begin{abstract}
We construct natural selfmaps of compact cohomgeneity one mani\-folds
with finite Weyl group and compute their degrees and Lefschetz numbers.
On manifolds with simple cohomology rings this yields in certain cases relations
between the order of the Weyl group and the Euler characteristic of a principal orbit.
We apply our construction to the compact Lie group $\mathrm{SU}(3)$
where we extend identity and transposition to an infinite family of selfmaps of
every odd degree. The compositions of these selfmaps with the power maps
realize all possible degrees of selfmaps of $\mathrm{SU}(3)$.
\end{abstract}

\maketitle

\section{Introduction}
A natural problem in topology is the following: Given a compact oriented mani\-fold $M^n$,
which integers can occur as the degree of maps $M\to M$? In the fundamental case
where $M$ is a sphere $\Sph^n$ all integers can easily be realized since $\Sph^n$ is a
$(n-1)$-fold suspension of $\Sph^1$. For other manifolds the problem is usually difficult
(see \cite{duan} for references and a detailed study in the case of
$(n-1)$-connected $2n$-manifolds).

In order to construct natural maps of degree $\neq -1,0,1$ one might first impose
actions of compact Lie groups $G$ on $M$ and then try to find equivariant maps.
In the most symmetric case where the action $G\times M \to M$ is transitive, it is
well-known and easy to see that every equivariant map is a diffeomorphism
and hence has degree~$\pm 1$. We deal with the case where the action is of
cohomogeneity one, i.e., the principal orbits $G/H$ are of codimension\,$1$
or, equivalently, the orbit space $M/G$ is $1$-dimensional and hence a
closed interval or a circle. We focus on the much more interesting case where the
orbit space $M/G$ is a closed interval.

We suppose that the cohomogeneity one manifold $M$ is equiped with an invariant
Riemannian metric such that the Weyl group is finite or, equivalently, such that the
normal geodesics are closed. Such metrics exist on many cohomogeneity one
manifolds $M$. For example, all invariant metrics with positive sectional curvature
necessarily have finite Weyl groups \cite{gwz}.
There might, however, be infinitely many invariant Riemannian metrics with
mutually distinct Weyl groups on the same compact manifold $M$ with a fixed
action $G\times M\to M$. The elementary but new construction of this paper
highly depends on the order of the Weyl group: We $k$-fold
the closed normal geodesics of $M$ starting from one of the non-principal orbits.
This leads to well-defined maps $\psi_k: M\to M$, $k= j\abs{W}/2+1$, for 
even integers $j$, and even for all integers $j$ provided a certain condition on the
isotropy groups is satisfied. In the following theorem we only deal with the general case.
For the exceptional case we refer to section\,\ref{degreeL}.

\begin{thOne}
A compact oriented Riemannian cohomogeneity one manifold $M$
whose orbit space $M/G$ is a closed interval and whose Weyl group $W$
is finite admits equivariant selfmaps $\psi_k : M\to M$
for all $k = j\abs{W}/2 +1$, $j\in 2\Z$ with~degree
\begin{gather*}
  \deg \psi_k = \begin{cases}
   k & \text{if the codimensions of both non-principal orbits are odd,}\\
  +1 & \text{otherwise,}
  \end{cases}
\end{gather*}
and Lefschetz number
\begin{gather*}
  L(\psi_k) = \begin{cases}
   - j \chi(G/H)/2 & \text{if the codimensions of both non-principal orbits are odd,}\\
  \phantom{-}\chi(M) & \text{otherwise.}
  \end{cases}
\end{gather*}
Here, $\chi(G/H)$ denotes the Euler characteristic of a principal orbit.
\end{thOne}

Degree and Lefschetz number of selfmaps are of course not independent.
In the case where $M^n$ is a rational homology sphere they are coupled
by the simple equation $L(\psi) = 1 + \deg\psi$ for even $n$ or $L(\psi) = 1 - \deg\psi$ for odd $n$.
Hence, we deduce from Theorem\,1 that $\rank G = \rank H$ and $\chi(G/H) = \abs{W}$
if the codimensions of both singular orbits are odd\footnote{Note that $\chi(G/H) = \abs{W}$
implies $W = N(H)/H$.
For the linear cohomogeneity one actions on spheres with $\rank G = \rank H$
the identity $W = N(H)/H$ was already observed in \cite{gwz}.}.
On manifolds $M$ with somewhat more complicated cohomology rings one
can still use Theorem\,1 to derive restrictions for cohomogeneity one actions
$G\times M\to M$ with $\rank G = \rank H$, i.e., $\chi(G/H) > 0$.
For such actions the dimensions of all orbits are necessarily even (hence the
dimension of $M$ is necessarily odd), and the Weyl groups are necessarily finite.
Using Dirichlet's theorem on prime numbers in arithmetic progressions
one obtains, for example, the following consequence of Theorem\,1:

\begin{crl}
Let $G\times M \to M$ be a cohomogeneity one action with $\rank G = \rank H$.
If $M$ has the rational cohomology ring of a product
\begin{gather*}
  \Sph^{l_1} \times \Sph^{l_2} \times \cdots \Sph^{l_m},\quad
  l_1 < l_2 < \ldots < l_m,
\end{gather*}
then $\chi(G/H) = 2^{m-1} \abs{W}$.
\end{crl}

In order to illustrate Theorem\,1 by examples other than spheres
one might scan the classifications \cite{neumann}, \cite{parker}, \cite{hoelscher}
for simply connected cohomogeneity one manifolds in dimensions $\le 7$
where the codimensions of both singular orbits are odd.
This scan reveals two nontrivial $\Sph^3$-bundles $M^7_1 \to \Sph^4$
and $M^7_2 \to \CP^2$ with sections. We briefly discuss these example
in section\,\ref{examples} and deduce from the refined version of Theorem\,1
that $M^7_1$ and $M^7_2$ admit equvariant selfmaps $\psi_k$
for all $k\in \Z$ with degree $k$ and Lefschetz number $2(1-k)$ and $3(1-k)$, respectively.
The space $M^7_1$ provides an example for Corollary\,2.

\smallskip

Our main example appears in dimension $8$: the compact Lie group $\SU(3)$.
Simple polynomial selfmaps of $\SU(3)$ are already given by
the power maps $\rho_k: A\mapsto A^k$. These maps are equivariant with respect to
conjugation and have degree $\deg \rho_k = k^2$ by a well-known result of Hopf \cite{hopf}.
Our construction yields additional selfmaps $\psi_k$ of $\SU(3)$ of every odd degree $k$
and Lefschetz number $0$ that extend identity and transposition to an infinite family.
These maps are equivariant with respect to the cohomogeneity one action
\begin{gather*}
  \SU(3) \times \SU(3) \to \SU(3), \quad (A,B) \mapsto ABA\tp.
\end{gather*}
Explicit formulas are given in section\,\ref{ex}. Combining the maps $\rho_k$
with the maps $\psi_k$ and a simple argument involving Steenrod squares
we prove
\begin{thTwo}
For any $m\in\N$ and $\ell \in \Z$ there is a selfmap of $\SU(3)$
with degree $4^m (2\ell + 1)$. For each of these selfmaps the entries
of the target matrix are real polynomials in the complex entries of the argument
matrix. Other integers do not appear as degrees of selfmaps of $\SU(3)$.
\end{thTwo}

\smallskip

The author would like to thank A.~Rigas for some useful discussions.

\bigskip

\section{Construction of the equivariant selfmaps}
Before we describe our construction we first summarize the necessary facts
about cohomogeneity one manifolds. We refer to the standard sources
\cite{mostert}, \cite{bredon}, \cite{alex} and the recent paper \cite{gwz}.

Let $M$ be a compact manifold on which a compact Lie group $G$ acts
with cohomogeneity one. Then $M/G$ is a circle or a closed interval.
In the former case all orbits are principal and $M$ is a $G/H$-bundle over $\Sph^1$
with structure group $N(H)/H$. In particular, $M$ has infinite fundamental group.
We only consider the case where $M/G$ is a closed interval.
The regular part $M_0$ of $M$ (i.e., the union of principal orbits) projects to the interior
of the interval and the two end points of the interval correspond to non-principal orbits
$N_0$ and $N_1$. The manifold $M$ is obtained by gluing the normal disk bundles
over $N_0$ and $N_1$ along their common boundary. In particular,
$\chi(M) = \chi(N_0) + \chi(N_1) - \chi(G/H)$.

Suppose that $M$ is equipped with an invariant Riemannian metric $\ip$.
After rescaling we can assume that $M/G$ is isometric to the unit interval $[0,1]$.
A normal geodesic is an (unparametrized) geodesic that passes through
all orbits perpendicularly. Through any point $p\in M$ there is at least one
normal geodesic and precisely one if $p$ is contained in a principal orbit.
The group $G$ acts transitively on the set of normal geodesics. We fix a
normal geodesic segment $\gamma:[0,1] \to M$ such that $\gamma(]0,1[)$ projects to
the interior of $M/G$ and $p_0 = \gamma(0) \in N_0$ and $p_1 = \gamma(1) \in N_1$.
This segment is a shortest curve from $N_0$ to $N_1$. It extends to a parametrized
normal geodesic $\gamma : \R \to M$. The isotropy groups $G_{\gamma(t)}$ are
constant along $\gamma$ except at the points $\gamma(t)$ with $t\in \Z$,
i.e., at points where $\gamma$ intersects the non-principal orbits.
We denote the generic isotropy group along $\gamma$ by $H$.

The Weyl group $W$ of $(M,\ip)$ is by definition the subgroup of elements of $G$
that leave $\gamma$ invariant modulo the subgroup of elements that fix
$\gamma$ pointwise. It is a dihedral subgroup of $N(H)/H$.
The geodesic $\gamma$ is periodic if and only if $W$ is finite.
More precisely, the order $\abs{W}$ of $W$ equals the number of times
that $\gamma$ passes through a fixed principal orbit before it closes.
To give an example of what might happen on a standard space consider
the diagonal $\SO(3)$-action on $\Sph^2\times \Sph^2(\sqrt{\alpha})$.
The normal geodesics are closed if and only if $\alpha$ is rational.
In this case, the Weyl group is the dihedral group $D_{\abs{p-q}}$
of order $2\abs{p-q}$ where $\alpha = p/q$ for positive integers
$p$ and $q$ with $\gcd(p,q) = 1$. Otherwise, the Weyl group is the
infinite dihedral group $D_{\infty}$. 

\smallskip

In the following we assume that the Weyl group of $(M,\ip)$ is finite, i.e.,
the fixed unit speed normal geodesic $\gamma$ and all their translated copies
$g \cdot \gamma$ are closed with period $\abs{W}$.

\smallskip

Our construction is based on the following elementary fact: Let
\begin{gather*}
  \psi_k : \Sph^1\to \Sph^1,\quad
  \lambda \mapsto \lambda^k
\end{gather*}
be the $k$-th power of $\Sph^1$.
Except in the trivial case where $k=1$ the fixed points of $\psi_k$ are precisely
the $\abs{k-1}$-roots of $1$ and $\psi_k$ has the property
$\psi_k(\lambda\cdot\lambda_0) = \psi_k(\lambda)\cdot \lambda_0$
for $\lambda_0$ with $\lambda_0^{k-1} = 1$.
In other words, if we have a circle of any length with a fixed base point $p_0$
then $k$-folding the intrinsic distance to $p_0$ is a well-defined map $\psi_k$
of the circle with $\abs{k-1}$ fixed points and $k$-folding the distance to
any of these fixed points leads to the same map $\psi_k$.

\begin{lem}
\label{wd}
The assignment $g\cdot \gamma(t) \mapsto g\cdot\gamma(kt)$ leads to a well-defined
smooth map $\psi_k :M\to M$ with $k = j\abs{W}/2+1$ for all $j\in 2\Z$, and even for all
$j\in\Z$ if the isotropy groups at $\gamma(t_0)$ and $\gamma(t_0+\abs{W}/2)$
for one and hence all odd integers $t_0$ are not just conjugate but equal.
\end{lem}

\begin{proof}
We first check that the assignment $\psi_k : g\cdot \gamma(t) \mapsto g\cdot\gamma(kt)$
yields a well-defined map on the unparametrized closed normal geodesic $\gamma(\R)$.
Suppose that $g\cdot \gamma(\R) = \gamma(\R)$.
Then $g\cdot \gamma(t) = \gamma(2t_0 \pm 2t)$ with $t_0\in \Z$.
A simple computation shows $\psi_k(g\cdot \gamma(t)) = g\cdot \psi_k(\gamma(t))$.
The geometric reasons for this equivariance of $\psi_k$ under the Weyl group are that
all the points $\gamma(2t_0)$ are fixed points of $\gamma(t)\mapsto \gamma(kt)$
(here the special form of $k$ is essential) and that the $k$-folding of a circle
is the same from every fixed point.
Now $\psi_k$ is well-defined on $\gamma(\R)$ and hence on the regular part $M_0$
of $M$ since through any point in $M_0$ there is only one (unparametrized) normal geodesic.
On $N_0$ the assignment $g\cdot \gamma(t)\mapsto g\cdot \gamma(kt)$
obviously leads to the identity. On $N_1$ we also get the identity provided
$j\in 2\Z$. If the isotropy groups at $\gamma(t_0)$ and $\gamma(t_0+\abs{W}/2)$
are equal for one odd integer $t_0$ then by the action of the Weyl group this
is true for all odd integers. If $j$ is an odd integer and $k = j\abs{W}/2+1$ then
\begin{gather*}
  g\cdot \gamma(t_0)\mapsto g\cdot\gamma(k t_0) = g\cdot \gamma(t_0+\abs{W}/2).
\end{gather*}
This actually defines an equivariant diffeomorphism $N_1 \to N_1$ or $N_1 \to N_0$
depending on whether $\abs{W}/2$ is even or odd. All in all we have shown that there
is a commutative diagram
\begin{gather*}
\begin{CD}
  \nu N_0 @>{\times k}>> \nu N_0\\
  @V{\exp}VV @VV{\exp}V\\
  M @>{\psi_k}>> M
\end{CD}
\end{gather*}
where $\nu N_0$ denotes the normal bundle of $N_0$. This implies the
smoothness of $\psi_k$.
\end{proof}

\bigskip

\section{Degree and Lefschetz number}
\label{degreeL}
Let $M$ be a compact Riemannian cohomogeneity one manifold with finite
Weyl group $W$ such that the orbit space $M/G$ is isometric to the interval $[0,1]$.
The manifold $M$ is orientable if and only if the principal orbits are orientable
and none of the non-principal orbits is exceptional and orientable (see \cite{mostert}).
In this section we assume that $M$ is orientable and equipped with a fixed orientation.

\begin{thm}
\label{Lefschetz}
The Lefschetz number of $\psi_k: M \to M$, $k=j\abs{W}/2+1$, is
\begin{gather*}
  L(\psi_k) = \begin{cases}
   - j\chi(G/H)/2 & \text{if $\codim N_0$ and $\codim N_1$ are both odd,}\\
  \phantom{-}\chi(M) & \text{otherwise,}
  \end{cases}
\end{gather*}
if $j$ is even and
\begin{gather*}
  L(\psi_k) = \begin{cases}
   - j \chi(G/H)/2 & \text{if $\rank G = \rank H$,}\\
   \phantom{-}\chi(N_0) & \text{otherwise,}
   \end{cases}
\end{gather*}
if $j$ is odd provided that $\psi_k$ is well-defined in this case.
\end{thm}

\begin{proof}
We perturb the map $\psi_k: M\to M$ to a map with only finitely many fixed points and
compute the Lefschetz number as the sum of fixed point indices.
For this pertubation we use the map $l_g : M\to M$, $p \mapsto g\cdot p$ where $g\in G$
is a general element, i.e., the closure of $\{g^m\,\vert\, m\in\Z\}$ is a maximal torus of $G$.
By a classical theorem of Hopf and Samelson \cite{samelson}, the restriction of $l_g$
to any of the orbits $G\cdot p$ has precisely $\chi(G/G_p)$ fixed points and each of these
fixed points has fixed point index one, i.e., $\det (\1 - A) > 0$ where $A$ is the derivative
of $l_g$ at the fixed point (note that our sign convention is different from that in \cite{samelson}).
We now consider the composition $l_g\circ\psi_k$ for a general element $g\in G$
sufficiently close to $1\in G$.
The equivariant map $\psi_k$ maps orbits to orbits. It is clear from the construction of
$\psi_k$ that only finitely many orbits $V_1,\ldots, V_m$ are mapped to themselves,
i.e., $\psi_k(V_i) = V_i$. On each of these orbits, $\psi_k$ is an equivariant diffeomorphism.
Hence, the restriction of $\psi_k$ to $V_i$ is either the identity or does not have fixed points.
The latter property is preserved under small perturbations of $\psi_k$.
Therefore, the map $l_g\circ\psi_k$ can only have fixed points in the orbits $V_i$
on which $\psi_k$ restricts to the identity, and the number of fixed points in $V_i$
is $\chi(V_i)$. Let $p$ be one of these fixed points in the orbit $V$. In order to compute
the fixed point index at $p$ we need the derivative of $l_g$ and $\psi_k$ at $p$.
Both derivatives clearly preserve both tangent and normal space to $V$ at $p$.
The derivative of $\psi_k$ is multiplication with $k$ on the normal space and
the identity on the tangent space. The derivative of $l_g$ is close to the identity
on the normal space and given by a matrix $A$ with $\det(\1 - A) > 0$ on the tangent space
by the above menitoned result of Hopf. Hence, if $B$ denotes the derivative of
$l_g \circ \psi_k$ at $p$ then the fixed point index at $p$ is given by
\begin{gather*}
  \det(\1 - B) = (\sgn(1-k))^{\codim V} = (-\sgn j)^{\codim V}.
\end{gather*}
The final step is now to sum these fixed point indices. We first consider the case
where $j$ is even. Then $\psi_k$ restricts to the identity on both $N_0$ and $N_1$.
Since $\psi_k$ has $\abs{k-1}$ fixed points on $\gamma([0,\abs{W}])$, we have
$\abs{k-1}/\abs{W} -1 = \abs{j}/2-1$ fixed points on $\gamma(]0,1[)$.
These fixed points correspond precisely to the orbits on which $\psi_k$
restricts to the identity. Hence we get
\begin{multline*}
  L(\psi_k) = L(l_g\circ\psi_k) = (-\sgn j)^{\codim N_0} \chi(N_0) + (-\sgn j)^{\codim N_0} \chi(N_0)\\
  +(\abs{j}/2-1) (-\sgn j)\, \chi(G/H).
\end{multline*}
Now suppose $j$ is odd. Then $\psi_k$ restricts to the identity on $N_0$ but not on $N_1$.
Since $\psi_k$ has $\abs{k-1}$ fixed points on $\gamma([0,\abs{W}])$, we have
$\abs{k-1}/\abs{W} -1/2 = (\abs{j}-1)/2$ fixed points on $\gamma(]0,1[)$. Hence,
\begin{gather*}
  L(\psi_k) = L(l_g\circ\psi_k) = (-\sgn j)^{\codim N_0} \chi(N_0)
  + \frac{\abs{j}-1}{2} (-\sgn j) \,\chi(G/H).
\end{gather*}
It is straightforward to simplify these two formulas for $L(\psi_k)$ using the facts that
the Euler characteristic of an odd dimensional manifold vanishes, that
$\chi(M) = \chi(N_0)  + \chi(N_1) - \chi(G/H)$, and that $\chi(G/H) > 0$ implies that
$\rank G = \rank H$ and hence that all orbits have even dimension.
In the case where the codimensions of both singular orbits are odd
one also has to use the fact that $\chi(N_0) = \chi(N_1) $ since both
fibrations $G/H \to N_i$ have fibers that are even dimensional spheres. 
\end{proof}

The degree of $\psi_k$ is the sum of the oriented preimages of a regular value.
Since the two non-principal orbits are mapped to themselves or one to the other
we just need to consider the regular part $M_0$ of $M$. In order to determine the
correct orientations we transfer the problem from $M_0$ to
$G/H\times (\R \smallsetminus\Z)$ by the map
\begin{gather*}
  \phi:G/H \times \R \to M,\quad
  (gH,t) \mapsto g\cdot \gamma(t).
\end{gather*}
The restriction of $\phi$ to each $G/H\times ]\ell,\ell+1[$ with $\ell \in \Z$
is a diffeomorphism. We equip $\R$ with the standard orientation and
$G/H$ with an orientation such the restriction of $G/H \times \,]0,1[\,\to M_0$ of $\phi$
is orientation preserving. This defines a standard product orientation on
$G/H\times \R$. For the rest of the paper, however, we will equip each
$G/H\times ]\ell,\ell+1[$ with the orientation inherited from $M_0$ by $\phi$.

\begin{lem}
\label{orient}
If $\codim N_0$ and $\codim N_1$ are both odd then all
$G/H\times \,]\ell,\ell+1[\,$ inherit the standard product orientation
from $M_0$ via $\phi$.

If $\codim N_0$ and $\codim N_1$ are both even, only the
pieces $G/H\times \,]2\ell,2\ell+1[\,$ inherit the standard product orientation
from $M_0$ via $\phi$.

If $\codim N_0$ is odd and $\codim N_1$ is even then only the
pieces $G/H\times \,]4\ell,4\ell+1[\,$ and $G/H\times \,]4\ell-1,4\ell[\,$
inherit the standard product orientation from $M_0$ via $\phi$.

If $\codim N_0$ is even and $\codim N_1$ is odd then only the
pieces $G/H\times \,]4\ell,4\ell+1[\,$ and $G/H\times \,]4\ell+1,4\ell+2[\,$
inherit the standard product orientation from $M_0$ via $\phi$.
\end{lem}

\begin{proof}
Let $\sigma_0$ denote the geodesic reflection along $N_0$, i.e.,
$\sigma_0$ maps each point $g\cdot \gamma(t)$ in $M\smallsetminus N_1$
to $g\cdot \gamma(-t)$. This equivariant diffeomorphism of $M\smallsetminus N_1$ corresponds to $-\id$ in the normal bundle of the orbit $N_0$ and hence preserves
the orientation of $M\smallsetminus N_1$ if and only if $\codim N_1$ is even.
We have the equivariant commutative diagram
\begin{gather*}
\begin{CD}
G/H\times \,]0,1[\, @>{(gH,t)\mapsto (gH,-t)}>> G/H \times \,]\!-1,0[ \\
@V{\phi}VV @V{\phi}VV \\
M_0 @>{\sigma_0}>> M_0
\end{CD}
\end{gather*}
Note that the map $\,]0,1[\,\to \,]\!-1,0[$, $t\mapsto -t$ reverses the orientation of $\R$.
Hence, $G/H\times \,]\!-1,0[\,$ inherits the same orientation from $M_0$ as
$G/H\times \,]0,1[\,$ if and only if $\sigma_0$ reverses orientation, i.e.,
if and only if $\codim N_0$ is odd.
\end{proof}

\begin{cor}
If $M$ is orientable and $\codim N_0$ and $\codim N_1$ do not have the same
parity, then the order $\abs{W}$ of the Weyl group $W$ of $M$ is divisible by $4$.
\end{cor}

\begin{thm}
\label{degree}
The map $\psi_k: M\to M$, $k= j\abs{W}/2 +1$ has degree
\begin{gather*}
  \deg \psi_k = \begin{cases}
  k & \text{if $\codim N_0$ and $\codim N_1$ are both odd,}\\
  +1 & \text{otherwise,}
  \end{cases}
\end{gather*}
if $j$ is even, and degree
\begin{gather*}
  \deg \psi_k = \begin{cases}
  k & \text{if $\codim N_0$ and $\codim N_1$ are both odd,}\\
  0 & \text{if $\codim N_0$ and $\codim N_1$ are both even, $\abs{W}\not\in 4\Z$,}\\
  -1 & \text{if $\codim N_0$ is even, $\codim N_1$ is odd, and $\abs{W}\not\in 8\Z$,}\\
  +1 & \text{otherwise,}
  \end{cases}
\end{gather*}
if $j$ is odd provided that $\psi_k$ is well-defined in this case (see Lemma\,\ref{wd}).
\end{thm}

\begin{proof}
The point $\gamma(\tau)$ with $\tau = 1/2$ is a regular value of the map $\psi_k$
whose $\abs{k}$ preimages are given by $\gamma(t_m)$ with
\begin{gather*}
  t_m = \frac{m \abs{W} + \tau}{k}, \text{ $m \in \Z$}.
\end{gather*}
In order to compute whether the differential of $\psi_k$ at $\gamma(t_m)$
preserves or reverses orientation we lift $\psi_k: M_0 \to M_0$ to the map
\begin{gather*}
  \tilde\psi_k: G/H \times \R \to G/H \times \R,\quad
  (gH,t) \to \bigl(gH,k t\bigr)
\end{gather*}
and answer the same question for the differential of $\tilde\psi_j$ at $(eH,t_m)$ by applying Lemma\,\ref{orient} (note again that we use the orientation on
$G/H\times (\R\smallsetminus\Z)$ induced from $M_0$ by $\phi$).
If $\codim N_0$ and $\codim N_1$ are both odd then all pieces $G/H\times ]\ell,\ell+1[$
inherit the standard product orientation. Hence, $\deg\psi_k = k$ and we are done.

Thus assume that $\codim N_0$ or $\codim N_1$ are even.
For each preimage $p$ of $\gamma(\tau)$ there is precisely one $m\in\Z$
such that $p = \gamma(t_m)$ and $t_m \in \,]0,\abs{W}[\,$.
We now have to count how many $t_m$ are contained in each
of the intervals $\,]\ell,\ell+1[\,$ for $\ell \in \{0,1,\ldots,\abs{W}-1\}$.
A simple computation shows that $t_m\in \,]\ell,\ell+1[\,$ is equivalent~to
\begin{gather*}
  \frac{j}{2}\ell + \frac{\ell - \tau}{\abs{W}} < m
    < \frac{j}{2}(\ell +1) + \frac{\ell +1 - \tau}{\abs{W}}
\end{gather*}
if $j \ge 0$, and equivalent to
\begin{gather*}
  \frac{j}{2}(\ell+1) + \frac{\ell+1 - \tau}{\abs{W}} < m
    < \frac{j}{2}\ell + \frac{\ell - \tau}{\abs{W}}
\end{gather*}
if $j < 0$. Suppose first that $j$ is even. For $\ell > 0$ we have
$\frac{\ell - \tau}{\abs{W}}\in \,]0,1[\,$ and $\frac{\ell +1 - \tau}{\abs{W}} \in \,]0,1[$.
Hence there are precisely $\abs{j}$ of the $t_m$ in each interval $\,]\ell,\ell+1[$
and $\abs{j+1}/2$ in the interval $\,]0,1[$. When counted with orientation
using Lemma\,\ref{orient} the sum of preimages of $\gamma(\tau)$ is hence $+1$
since there are $\abs{W}$ intervals $]\ell,\ell+1[$ and $\abs{W}$ is
divisible by $2$ if the codimensions of $N_0$ and $N_1$ are both even
and $\abs{W}$ is divisible by $4$ if the codimensions of $N_0$ and $N_1$
have different parities. Note that in the case $j \ge 0$ there is one point
more in the interval $\,]0,1[$ than in the other intervals and this point has
positive orientation, while in the case $j<0$ there is one point less but the
orientation of the real line has changed. Suppose now that $j$ is odd.
Then $\frac{\ell - \tau}{\abs{W}}\in \,]0,\frac{1}{2}[\,$ if and only if $0 < \ell \le \abs{W}/{2}$
and $\frac{\ell +1 - \tau}{\abs{W}} \in \,]0,\frac{1}{2}[$ if and only if $0 \le \ell < \abs{W}/{2}$.
Using these facts it is straightforward to show that there are
$\abs{j+1}/2$ of the $t_m$ in each of the intervals $\,]\ell,\ell+1[\,$
with $\ell =0$, $\ell = \abs{W}/2$, odd $l < \abs{W}/2$, or even
$l > \abs{W}/2$, while there are $\abs{j-1}/2$ of the $t_m$ in each of the intervals
$]\ell,\ell+1[$ with even $0 < l < \abs{W}/2$ or odd $l > \abs{W}/2$.
The claimed degrees follow now by applying Lemma\,\ref{orient}.
\end{proof}

\begin{proof}[Proof of Corollary\,2]
Suppose that $M$ has the rational cohomology ring of
\begin{gather*}
  \Sph^{l_1} \times \Sph^{l_2} \times \cdots \Sph^{l_m},\quad
  l_1 < l_2 < \ldots < l_m,
\end{gather*}
and that $G\times M\to M$ is a cohomogeneity one action with $\rank G = \rank H$,
i.e., $\chi(G/H) > 0$ where $H$ is a principal isotropy group. By Lemma\,\ref{wd},
$M$ admits equivariant selfmaps $\psi_k$, $k=j\abs{W}/2+1$, with degree $k$
and Lefschetz number $L(\psi_k) = -j\chi(G/H)$ for every even integer $j$.
By Dirichlet's theorem there are infinitely many prime numbers in the arithmetic
progression $k=j\abs{W}/2+1$. Hence, we can assume that $k$ is a prime number.
We now inspect the induced map $\psi_k^{\ast}$ on the cohomology ring.
Let $x_1,\ldots, x_m$ be generators of $H^{\ast}(M)$ with $x_i \in H^{l_i}(M)$.
The equation $\deg\psi_k = k$ means $\psi_k^{\ast}(x_1x_2\cdots x_m) = k x_1x_2\cdots x_m$.
Since $k$ is a prime number and $l_1 < l_2 < \ldots < l_m$ it is clear that
\begin{gather*}
  \psi_k^{\ast}(x_j) = k \sigma_j x_j + \text{products of lower dimensional generators}
\end{gather*}
with $\sigma_j = \pm 1$ for precisely one $j \in \{1,\ldots,m\}$ and
\begin{gather*}
  \psi_k^{\ast}(x_i) = \sigma_i x_i + \text{products of lower dimensional generators}
\end{gather*}
with $\sigma_i = \pm 1$ for $i \in \{1,\ldots,m\}$, $i\neq j$. The total number of $\sigma_i$ with
$\sigma_i = -1$ for $i  \in \{1,\ldots,m\}$ must be even.
The Lefschetz number $L(\psi_k)$ is equal to $-j\chi(G/H)/2$ by Theorem\,\ref{Lefschetz}.
On the other hand, $L(\psi_k)$ is the alternating trace of $\psi^{\ast}$ in cohomology.
The simple cohomology ring structure implies
\begin{gather*}
  -j\chi(G/H)/2 = L(\psi_k) = \bigl(1+(-1)^{l_j} k\sigma_j \bigr)
    \prod_{i\neq j} \bigl(1+(-1)^{l_i} \sigma_i \bigr).
\end{gather*}
Since we can assume $j > 2$ and since $\chi(G/H) > 0$, we have
$\sigma_i = (-1)^{l_i}$ for all $i\neq j$ and $\sigma_j = -(-1)^{l_j}$.
Hence, $\chi(G/H) = 2^{m-1}\abs{W}$ and Corollary\,2 is proved.
\end{proof}

\bigskip

\section{The main example $\SU(3)$}
\label{ex}
The natural selfmaps of compact Lie groups $G$ are the power maps
$\rho_k:A\mapsto A^k$. These maps are equivariant with respect to conjugation
and hence determined by their restriction to a maximal torus
$T^r\subset G$. Using this fact it is easy to see that the degree of $\rho_k$
is $k^r$ and thus rather few integers can be realized as degrees of
power maps. For $\SU(3)$ with rank $r = 2$ these are all squares $k^2$.

Conjugation on $\SU(3)$ is an action of cohomogeneity two. We use
the cohomogeneity one action
\begin{gather*}
  \SU(3) \times \SU(3) \to \SU(3),\quad (A,B) \mapsto A B A\tp
\end{gather*}
to construct polynomial selfmaps of $\SU(3)$ whose degrees cover all odd integers
and that extend identity and transposition to an infinite family.
Polynomial here means that the entries of the target matrix are real polynomials
in the complex entries of the argument matrix.

\begin{lem}
$\SU(3)$ admits selfmaps $\psi_k$ with degree $k$ for all odd integers $k$.
These maps are equivariant with respect to the cohomogeneity one action above.
\end{lem}

\begin{proof}
A straightforward computation shows that 
\begin{gather*}
  \gamma(t) = \bmat \cos t & -\sin t & \m 0\\ \sin t & \m \cos t & \m 0\\
    0 & \m 0 & \m 1\emat
\end{gather*}
is a geodesic that passes perpendicularly through all orbits.
This geodesic is closed with period $2\pi$. At times $t = 0$ and $t=\pi$ it passes
through the singular orbit $\SU(3)/\SO(3)$ that consists of the symmetric matrices
in $\SU(3)$. At $t = \frac{\pi}{2}+\pi\Z$ it passes through the other singular orbit
$\SU(3)/\SU(2)$ where $\SU(2)$ is embedded in the upper left corner of $\SU(3)$.
The isotropy group of $\gamma(t)$ for all other times $t$ is $\SO(2)$ embedded
in the upper left corner of $\SU(3)$ in the standard way.
Note that the singular isotropy groups at $t=\frac{\pi}{2}$ and $t=\frac{3\pi}{2}$
are not just conjugate but equal.
The claim follows from Lemma\,\ref{wd} and Theorem\,\ref{degree}.
\end{proof}

We now derive an explicit formula for these maps. Let
\begin{align*}
  f_j(x) &= \sum_{i=0}^j (-1)^i \binom{2j+1}{2i} x^i (1-x)^{j-i},\\
  g_j(x) &= \sum_{i=0}^j (-1)^i \binom{2j+1}{2i+1} x^i (1-x)^{j-i},\\
  h_j(x) &= \frac{1-f_j(x)}{x} =
    \sum_{i=0}^{j-1} (-1)^i \biggl( \binom{j}{i+1}x^i +  \binom{2j+1}{2i+2} x^i (1-x)^{j-1-i} \biggr).
\end{align*}

\begin{thm}
The selfmaps $\psi_k :\SU(3) \to \SU(3)$ of odd degree $k=2j+1$ are given by the
following formula:
\begin{gather*}
  \psi_{2j+1}(B) = \tfrac{1}{2} f_{\abs{j}}(x)(B+B\tp) \pm \tfrac{1}{2} g_{\abs{j}}(x) (B-B\tp)
  + \tfrac{1}{4} h_{\abs{j}}(x) \bmat \bar b_{23} - \bar b_{32}\\ \bar b_{31} - \bar b_{13}\\
  \bar b_{12} - \bar b_{21} \emat \,
  \bmat \bar b_{23} - \bar b_{32}\\ \bar b_{31} - \bar b_{13}\\
  \bar b_{12} - \bar b_{21} \emat\tp
\end{gather*}
where $x = 1 - \tfrac{1}{4}(\tr B - 1)^2$. The $+$-sign has to be taken
if $j \ge 0$ and the $-$-sign has to be taken if $j < 0$.
\end{thm}

\begin{proof}
Let $B = (b_{ij})$ be a fixed matrix in $\SU(3)$. We want to solve the equation
$B = A\gamma(t)A\tp$ for $t\in [0,\tfrac{\pi}{2}]$ and $A\in \SU(3)$ then
multiply $t$ by an odd integer $k$ and determine $B' = A\gamma(kt) A\tp$.
Symbolically, we illustrate this as follows
\begin{gather*}
  B \rightsquigarrow (A,t) \rightsquigarrow (A,kt) \rightsquigarrow B'.
\end{gather*}
We will see that it is not necessary to determine the entries of $A$ completely with all
their ambiguity. We stop at an intermediate level and then return.

Let $A = (z,w,v)$ with $z,w,v\in \C^3$. Then $B = A\gamma(t)A\tp$ becomes
\begin{gather*}
  b_{ii} = (z_i^2 + w_i^2)\cos t + v_i^2,\\
  b_{i,i+1} = (z_iz_{i+1}+w_iw_{i+1}) \cos t + v_i v_{i+1} - \bar v_{i-1}\sin t,\\
  b_{i,i-1} = (z_iz_{i-1}+w_iw_{i-1}) \cos t + v_i v_{i-1} + \bar v_{i+1}\sin t\\
\end{gather*}
with indices cyclic modulo 3. It is straightforward to see that
\begin{gather*}
  \cos t = \tfrac{1}{2} ( \tr B -1),\\
  2\bar v_1 \sin t = b_{32}-b_{23},\\
  2\bar v_2 \sin t = b_{13}-b_{31},\\
  2\bar v_3 \sin t = b_{21}-b_{12},\\
  (z_i^2 + w_i^2) \cos t + v_i^2 = b_{ii},\\
  (z_2z_3+w_2w_3) \cos t + v_2 v_3 = \tfrac{1}{2} (b_{23}+b_{32}),\\
  (z_3z_1+w_3w_1) \cos t + v_3 v_1 = \tfrac{1}{2} (b_{31}+b_{13}),\\
  (z_1z_2+w_1w_2) \cos t + v_1 v_2 = \tfrac{1}{2} (b_{12}+b_{21}).
\end{gather*}
Now, applying the binomial formula to $(\cos t + \cplxi \sin t)^k$
leads the following formulas for $\cos kt$ and $\sin kt$ with $k$ odd:
\begin{gather*}
  \cos (2j+1)t = f_j(\sin^2 t) \cos t, \qquad \sin (2j+1)t = g_j(\sin^2 t) \sin t.
\end{gather*}
We obtain
\begin{align*}
  b_{11}' &= (z_1^2 + w_1^2) \cos (2j+1)t + v_1^2\\
  &=  (z_1^2 + w_1^2) f_j(\sin^2 t) \cos t + v_1^2 f_j(\sin^2 t) + v_1^2 h_j(\sin^2 t) \sin^2 t\\
  &= b_{11} f_j(\sin^2 t) + \tfrac{1}{4} (\bar b_{23} - \bar b_{32})^2 h_j(\sin^2 t)
\end{align*}
and analogous expressions for $b_{22}$ and $b_{33}$. Note that these are polynomials
in the matrix entries of $B$ since $\sin^2 t = 1 - \tfrac{1}{4}(\tr B - 1)^2$ is a polynomial in
$\tr B$ and, e.\,g., $\bar b_{23} = b_{31}b_{12} - b_{11}b_{32}$, since $B\in \SU(3)$.
Similarly, we obtain polynomial expressions for the off-diagonal terms of
$B' = \psi_{2j+1}(B)$. All in all this leads to the claimed formula.
\end{proof}

\begin{cor}
For every $m\in\N$ and every odd integer $(2\ell +1)$ there is a
polynomial selfmap $\SU(3)\to \SU(3)$ of degree $4^m (2\ell +1)$.
\end{cor}

\begin{proof}
The maps $\psi_k$ are defined for odd $k$ and have degree $k$.
The power maps $\rho_k: A\mapsto A^k$ have degree $k^2$.
\end{proof}

\begin{lem}
Let $\phi: \SU(3)\to \SU(3)$ be any map.
Then $\deg \phi = 4^m \cdot (\text{odd number})$ for some $m\in\N$.
\end{lem}

\begin{proof}
The cohomology ring $H^{\ast}(\SU(3); \Z)$ is isomorphic to the
cohomology ring $H^{\ast}(\Sph^3\times\Sph^5; \Z)$, i.e., to
$\Lambda(x,y)$ with generators $x$ in dimension~$3$ and $y$ in dimension~$5$.
If $\phi: \SU(3) \to \SU(3)$ is any map, we have the commutative diagram
\begin{gather*}
\begin{CD}
  H^3(\SU(3); \Z_2) @>{\phi^{\ast}}>> H^3(\SU(3); \Z_2)\\
  @V{\Sq^2}VV  @V{\Sq^2}VV\\
  H^5(\SU(3); \Z_2) @>{\phi^{\ast}}>> H^5(\SU(3); \Z_2)
\end{CD}
\end{gather*}
Vertically, the second Steenrod square $\Sq^2$ yields an isomorphism
(see \cite{bredtop}, p.\,469).
The diagram thus says that $\phi^{\ast}(y)$ is an even multiple of $y$ if and only if
$\phi^{\ast}(x)$ is an even multiple of $x$. Hence, $\phi^{\ast}(xy)$ is in
$4^m (2\ell + 1)\cdot xy$ for some $m\in \N$ and $\ell \in \Z$.
\end{proof}

\begin{lem}
We have $\rho_k^{\ast}(x) = k x$ and $\rho_k^{\ast}(y) = k y$ while
$\psi_k^{\ast}(x) = k x$ and $\psi_k^{\ast}(y) = y$. Hence, none of the
$\rho_k$ is homotopic to any of the $\psi_k$.
\end{lem}
\begin{proof}
For $\rho_k$ we refer to \cite{hopf}. It is easy to see that the subgroup $\SU(2)$
embedded in the upper left corner of $\SU(3)$ is invariant under $\psi_k$
and $\psi_k$ restricts to a selfmap of $\SU(2)$ with degree $k$.
The commutative diagram
\begin{gather*}
\begin{CD}
H^3(\SU(3);\Z) @>{\psi_k}>> H^3(\SU(3);\Z)\\
@VVV @VVV\\
H^3(\SU(2);\Z) @>{\times k}>> H^3(\SU(2);\Z)
\end{CD}
\end{gather*}
shows $\psi_k^{\ast}(x) = kx$. Since $\deg \psi_k = k$, the generator
$xy$ of $H^8(\SU(3;\Z)$ is mapped to $k xy$. Hence, $\psi_k^{\ast}(y) = y$.
\end{proof}

\bigskip

\section{Further examples}
\label{examples}
In this section we supply examples for the various cases in Theorem\,\ref{Lefschetz}
and Theorem\,\ref{degree}.
Classical examples of selfmaps are the $k$-powers of $\Sph^n\subset \R\times \R^n$,
\begin{gather*}
  (\cos t,v \sin t) \mapsto (\cos kt, v\sin kt).
\end{gather*}
These well-defined analytic selfmaps
of $\Sph^n$ are equivariant with respect to the standard $\OO(n)$-action on $\R^n$.
If $n$ is odd we have $\rank \OO(n) = \rank \OO(n-1)$. From Theorem\,\ref{degree}
we recover the known fact that the degree of the $k$-power of $\Sph^n$ is $k$ if
$n$ is odd, and $0$ or $1$ if $n$ is even, depending on whether $k$ is even or odd.
From Theorem\,\ref{Lefschetz} we recover the Lefschetz number $(1-k)$ if $n$ is odd,
and $1$ or $2$ if $n$ is even, depending on whether $k$ is even or odd.

An example where the codimensions of the singular orbits have different parity
and there are equivariant maps with degree $-1$ is given by the standard action of
$\SO(1+m)\subset \SU(1+m)$ on $\CP^m$ for odd $m$.
A normal geodesic is given by
\begin{gather*}
  \gamma(t) = [\cos t : \cplxi \sin t : 0 : \ldots : 0].
\end{gather*}
The singular orbit at $t=0$ is $\RP^m = \SO(1+m)/\OO(m)$ and the singular
isotropy groups at $t = \pi/4$ and $t=-\pi/4$ are both equal to $\SO(2)\SO(m-1)$.
The normal geodesic $\gamma$ is closed with period $\pi$.
Hence, the Weyl group is the dihedral group of order~$4$ and $(2j+1)$-folding the
distance to $N_0 = \RP^m$ yields a well-defined map with degree $1$ and Lefschetz
number $m+1$ if $j$ is even and degree $-1$ and Lefschetz number $0$ if $j$ is odd.
Note that for even $m$ selfmaps of $\CP^m$ with degree $-1$ cannot exist and this
is perfectly matched by the fact that both singular orbits have even codimensions in this case.

\smallskip

We finally discuss the spaces $M^7_1 = \Syp(2)\times_{\Syp(1)^2}\Syp(1)$
and $M^7_2 = \SU(3)\times_{\U(2)}\SU(2)$ where $\U(2)$ acts on $\SU(2)$
by conjugation and one of the factors of $\Syp(1)^2$ acts trivially on $\Syp(1)$
and the other by conjugation. The space $M^7_1$ is an $\Sph^3$-bundle over $\Sph^4$.
Such bundles are classified by their characteristic homomorphism $\Sph^3\to \SO(4)$.
Using two Cartan embeddings of $\Sph^4$ into $\Syp(2)$ it is not difficult to see that
for $M^7_1$ this homomorphism is given by $q \mapsto C_q$ where $q$ is a unit
quaternion and $C_q: \H \to \H$ is conjugation by $q$. In particular, the bundle
$M^7_1\to \Sph^4$ has a section but is nontrivial. It is known that $M^7_1$ is not
homotopy equivalent to $\Sph^3\times \Sph^4$ but has the same cohomology and
homotopy groups as $\Sph^3\times \Sph^4$ \cite{james}.
The group $\Syp(2)$ acts from the left on $M^7_1$ by cohomogeneity one.
The principal orbits are diffeomorphic to $\CP^3 = \Syp(2)/\Syp(1)\SO(2)$
with Euler characteristic $4$ and the singular orbits are both diffeomorphic to $\Sph^4$.
We have $\rank G = \rank H$. The Weyl group is isomorphic to $\Z_2$.
Hence, by Lemma\,\ref{wd} there exist equivariant selfmaps $\psi_k$ with degree
$\psi_k = k$ and Lefschetz number $L(\psi_k) = 2(1-k)$ for all integers $k$.
This example thus illustrates Corollary\,2, since $4 = \chi(G/H) = 2\cdot \abs{W}$.

The space $M^7_2$ is an $\Sph^3$-bundle over $\CP^2$.
The group $\SU(3)$ acts fromt the left on $M^7_2$ by cohomogeneity one.
The principal orbits are diffeomorphic to $\SU(3)/\T^2$ with Euler characteristic $6$
and the singular orbits are both diffeomorphic to $\CP^2$. Since the bundle
$M^7_2\to \CP^2$ has a section, $M^7_2$ has the same cohomology and homotopy groups
as $\Sph^3 \times \CP^2$. Moreover, the corresponding principal $\SO(3)$-bundle over
$\CP^2$ is the Aloff-Wallach space $W^7_{1,1}$. The first Ponrjagin class of this
principal bundle is $-3\in H^4(\CP^2)\approx \Z$ \cite{ziller}. Hence,
the $\Sph^3$-bundle $M^7_2 \to \CP^2$ has $w_2 = 0$ and $p_1 = -3$.
By the Dold-Whitney classification \cite{dold} it is nontrivial.
The Weyl group of the cohomogeneity one $\SU(3)$ action is isomorphic to $\Z_2$.
Hence, there exist equivariant selfmaps $\psi_k$ with degree $\psi_k = k$ and
Lefschetz number $L(\psi_k) = 3(1-k)$ for all integers $k$ in agreement with the
homology of $M^7_2$.

\bigskip

%
%


\nocite{*}
\begin{thebibliography}{DoWh}
%
%
\bibitem[AA]{alex}
A.\,V.~Alekseevsky, D.\,V.~Alekseevsky, {\em Riemannian $G$-manifolds with
one-dimensional orbit space}, Ann.\ Global Anal.\ Geom. {\bf 11} (1993),
197--211.
%
\bibitem[Bd1]{bredon}
G.\,E.~Bredon, {\em Introduction to compact transformation groups},
Pure and Applied Mathematics, Vol. 46, Academic Press, New York, 1972.
%
\bibitem[Bd2]{bredtop}
G.\,E.~Bredon, {Topology and Geometry}, GTM 139, Springer, New York 1993.
%
\bibitem[DoWh]{dold}
A.~Dold, H.~Whitney, {Classification of oriented sphere bundles over a 4-complex},
Ann.\ Math. {\bf 69} (1959), 667-677.
%
\bibitem[DuWa]{duan}
H.~Duan, S.~Wang, {\em The degree of maps  between manifolds},
Math.\ Z. {\bf 244} (2003), 67--89.
%
\bibitem[GWZ]{gwz}
K.~Grove, B.~Wilking, W.~Ziller, {\em Positively curved cohomogeneity one manifolds and $3$-Sasakian geometry}, J.\ Different.\ Geom., to appear.
%
\bibitem[Hoe]{hoelscher}
C.~Hoelscher, {\em Classification of cohomogeneity one manifolds in low dimensions},
Ph.\ D. thesis, University of Pennsylvania, Philadelphia 2007.
%
\bibitem[Ho]{hopf}
H.~Hopf, {\em \"Uber den Rang geschlossener Liescher Gruppen},
Comment.\ Math.\ Helv. {\bf 13} (1941), 119--143.
%
\bibitem[HS]{samelson}
H.~Hopf, H.~Samelson, {\em Ein Satz \"uber die Wirkungsr\"aume geschlossener Liescher Gruppen},
Comment.\ Math.\ Helv. {\bf 13} (1941), 240--251.
%
\bibitem[JW]{james}
I.\,M.~James, J.\,H.\,C.~Whitehead, {\em The homotopy theory of sphere bundles over spheres.~I.},
Proc.\ Lond.\ Math.\ Soc. {\bf 4} (1954), 196--218.
%
\bibitem[Mo]{mostert}
P.\,S.~Mostert, {\em On a compact Lie group acting on a manifold},
Ann.\ of Math. {\bf 65} (1957), 447-455; Errata ibid. {\bf 66}, 589.
%
\bibitem[Neu]{neumann}
W.\,D.~Neumann, ``$3$-dimensional $G$-manifolds with $2$-dimensional orbits'',
pp. 220--222 in {\em Proc. Conf. on Transformation Groups}
(New Orleans, La., 1967), Springer, New York, 1968
%
\bibitem[Pa]{parker}
J.~Parker, {\em $4$-dimensional $G$-manifolds with $3$-dimensional orbits},
Pac.\ J.\ Math. {\bf 125} (1986), 187-204.
%
\bibitem[Zi]{ziller}
W.~Ziller, {\em Fatness revisited}, unpublished lecture notes, University of Pennsylvania, 2000.
\end{thebibliography}
\end{document}